\documentclass[12pt]{article}
\usepackage{amsmath}
\usepackage{amssymb}
\usepackage{dsfont}
\usepackage{cite}
\newsymbol \blackbox 1004
\newcommand{\eh}{\hfill}\newlength{\sperr}

\newenvironment{proof}{{\settowidth{\sperr}{\bf\rm
Proof}%
\par\addvspace{0.3cm}\noindent\parbox[t]{1.3\sperr}
{\bf\rm P\eh r\eh o\eh o\eh f\eh }%
}}{\nopagebreak\mbox{}
$\blackbox$\par\addvspace{0.3cm}}

\def\nn{\nonumber}

\def\vk{\varkappa}

\def\s{\sigma}
\def\la{\lambda}

\def\vp{\varphi}
\def\vt{\vartheta}
\def\ve{\varepsilon}
\def\wh{\widehat}
\def\wt{\widetilde}
\def\ov{\overline}

\def\BC{{\mathbb C}}
\def\BR{{\mathbb R}}
\def\BN{{\mathbb N}}
\def\clp{{\mathcal P}}
\def\cla{{\mathcal A}}
\def\clb{{\mathcal B}}
\def\clc{{\mathcal C}}
\def\cld{{\mathcal D}}
\def\clg{{\mathcal G}}

\def\cln{{\mathcal N}}

\def\cld{{\mathcal D}}

\def\spa{{\rm Span}}

\newcommand{\I}{\mathrm{i}}

\newtheorem{Pa}{Paper}[section]
\newtheorem{Tm}[Pa]{{\bf Theorem}}

\newtheorem{Rk}[Pa]{{\bf Remark}}

\newtheorem{Dn}[Pa]{{\bf Definition}}

\newtheorem{Pn}[Pa]{{\bf Proposition}}

\newenvironment{dedication}
        {\vspace{1ex}\begin{quotation}\begin{center}\begin{em}}    
        {\par\end{em}\end{center}\end{quotation}}

\title{General-type discrete self-adjoint  Dirac systems:  explicit solutions of direct and inverse
problems, \\ asymptotics of Verblunsky-type coefficients and stability of solving inverse problem.}

\author{I.Ya. Roitberg and A.L. Sakhnovich}

\date{}
\begin{document}
\maketitle

\begin{dedication} To V.A. Marchenko with admiration
\end{dedication}

\begin{abstract}    We consider discrete self-adjoint  Dirac systems determined by the
potentials (sequences) $\{C_k\}$ such that the matrices $C_k$ are positive definite
and $j$-unitary, where $j$ is a  diagonal $m\times m$ matrix and has $m_1$ entries $1$ and $m_2$ entries
$-1$ ($m_1+m_2=m$) on the main diagonal. We construct systems with rational Weyl functions and explicitly
solve inverse problem to recover systems from the contractive rational Weyl functions.
Moreover, we study the stability of this procedure. The matrices $C_k$ (in the potentials) are so called
Halmos extensions of the Verblunsky-type coefficients $\rho_k$. We show that in the case of 
the contractive rational Weyl functions the coefficients $\rho_k$ tend to zero and the matrices $C_k$ tend
to the indentity matrix $I_m$.
\end{abstract}


\vspace{0.2em}


\section{Introduction}\label{Intro}
\setcounter{equation}{0}
General-type 
discrete self-adjoint Dirac systems have the form:
\begin{equation} \label{0.1}
y_{k+1}(z)=(I_m+ \I z
j C_k)
y_k(z) \quad \left( k \in \BN_0
\right),
\end{equation}
where $\BN_0$ stands for the set of non-negative integers,  $I_m$ is the $m \times m$ identity matrix, $"\I"$ is the imaginary unit
($\I^2=-1$) and the $m \times m$ matrices $\{C_k\}$ are positive and  $j$-unitary:
\begin{equation} \label{0.2}
C_k>0, \quad C_k j C_k=j, \quad  j: = \left[
\begin{array}{cc}
I_{m_1} & 0 \\ 0 & -I_{m_2}
\end{array}
\right] \quad (m_1+m_2=m; \, \, m_1, \, m_2 \not= 0).
\end{equation}
First, we will consider (in Section \ref{GBDT})  explicit solutions of the direct and inverse problems for system \eqref{0.1}, \eqref{0.2}
in terms of Weyl-Titchmarsh (or simply Weyl) functions. General-type direct and inverse problems 
for this system were studied (in terms of Weyl functions) in \cite{FKRS14}
and explicit solutions in the case $m_1=m_2$ were dealt with in \cite{FKRS08}. Our Section \ref{GBDT} (and Appendix) 
complete the results from \cite{FKRS14} by adding the properties  of the Weyl functions in the lower
half-plane and  generalize the explicit results from \cite{FKRS08} for the case when $m_1$ does not necessarily equal $m_2$.
We will often reduce our proofs in Section \ref{GBDT} and Appendix and  refer
to the more detailed proofs in \cite{FKRS08, FKRS14}. However, a complete procedure of explicitly solving the inverse
problem from Section \ref{GBDT} is missing in \cite{FKRS08} (and so it is new for $m_1=m_2$ as well).

The case of explicit solutions of direct and inverse problems corresponds to the {\it rational Weyl functions}.
The results in Section \ref{GBDT} are based on our generalized B\"acklund-Darboux (GBDT) approach,
which was initiated by the seminal book \cite{Mar} by V.A. Marchenko. For various versions of B\"acklund-Darboux transformations and related
commutation methods see, for instance, \cite{Ci, D, GeT, KoSaTe, MS, KS, ALS94, SaSaR} and references therein.

Section \ref{Verb} is dedicated to the asymptotics of the {\it potentials} (sequences) $\{C_k\}$ corresponding to rational
Weyl functions. For this purpose, we first derive the asymptotics of the so called \cite{ALSverb} Verblunsky-type coefficients.

Finally, in Section \ref{stab}, we study stability of  our method of explicit solving inverse problem
for system \eqref{0.1}, \eqref{0.2},
and these results are new even in the cases $m_1=m_2$ and $m_1=m_2=1$.
We note that various important early results on the stability of solving  inverse problems
were obtained by V.A.~Marchenko (see, e.g., \cite{Mar0}).

 In the paper, $\BN$ denotes the set of natural numbers, $\BR$ denotes the real axis, $\BC$ stands for the complex plane, and
$\BC_+$ ($\BC_-$) stands for the open upper (lower) half-plane. The spectrum of a square matrix $A$ is denoted by $\s(A)$.
\section{GBDT and direct and inverse problems}\label{GBDT}
\setcounter{equation}{0}
\paragraph{1.} The fundamental $m \times m$
solution $\{W_k \}$ of  \eqref{0.1}  is normalized by
\begin{align} \label{0.3}&
W_0(z)=I_m.
\end{align}
For the case $z \in \BC_+$, the definition of the Weyl function $\vp(z)$  of Dirac system \eqref{0.1}, \eqref{0.2} 
was given in  \cite{FKRS14} in terms of $W_k(z)$. Below we define the Weyl function in $\BC_-$, which is somewhat more convenient for 
our purposes. Clearly, this Weyl function has similar properties to those in  \cite[Theorem 3.8]{FKRS14}.
\begin{Dn} \label{defWeyl}  The Weyl function of the Dirac system \eqref{0.1} 
$($which is given on the semi-axis $0\leq k < \infty$ and satisfies \eqref{0.2}$)$ is an $m_1\times m_2$ matrix  function $\vp(z)$ in the
lower half-plane, such that the following inequalities hold:
\begin{align} \label{0.4}&
\sum_{k=0}^\infty q(z)^k 
\begin{bmatrix}
 \vp(z)^* &I_{m_2} 
\end{bmatrix}
W_k(z)^*C_k W_k (z)
\begin{bmatrix}
 \vp(z) \\ I_{m_2} 
\end{bmatrix}<\infty \quad (z\in \BC_-),
\\ & \label{0.5}
q(z):=(1+|z|^2)^{-1}.
\end{align}
\end{Dn}
The properties of the Weyl function are described in the theorem below, which is proved in Appendix (using the standard
Weyl disk procedure).
\begin{Tm} \label{TmWf}
There is a unique Weyl function of  the discrete Dirac system \eqref{0.1}, 
which  is given on the semi-axis $ 0\leq k < \infty$ and satisfies \eqref{0.2}. 
This Weyl function $\vp$ is analytic and contractive $($i.e.,
$\vp^*\vp \leq I_{m_2})$ on $\BC_-$.
\end{Tm}
In the proof of Theorem \ref{TmWf} in Appendix, we will  need the inequalities
\begin{equation} \label{Prop2.2}
C_k\geq j,
\end{equation}
which (together with the inequalities $\, C_k \geq -j\,$) immediately follow from \cite[Proposition 2.2]{FKRS14}.

Another way to prove Theorem \ref{TmWf} and the uniqueness of the solution of the inverse problem, which we will need
further, is to consider  Dirac systems
\begin{align} \label{0.1'} &
\wt y_{k+1}(z)=(I_m+ \I z \,
\wt j \, \wt C_k)
\wt y_k(z) \quad \left( k \in \BN_0
\right),
\\ \label{d1} &
\wt j:=-JjJ^*=\begin{bmatrix}I_{m_2} & 0\\ 0 & -I_{m_1} \end{bmatrix}, \quad J:=\begin{bmatrix}0 & I_{m_2} \\ I_{m_1} & 0\end{bmatrix},
\quad \wt C_k:=JC_k J^*.
\end{align}
Systems \eqref{0.1'}, \eqref{d1} are dual to the systems \eqref{0.1}, \eqref{0.2}, and it is immediate from \eqref{0.2}, \eqref{d1}
that the relations
\begin{equation} \label{0.2'}
J^*J=I_m, \quad \wt C_k>0, \quad \wt C_k \, \wt  j \,\wt C_k=\wt j
\end{equation}
are valid. Hence, systems \eqref{0.1'} are again self-adjoint Dirac systems. Similar to $\wt j$ and $\wt C_k$
we use ``tilde" in other notations (introduced for self-adjoint Dirac systems), when it goes about
systems \eqref{0.1'}. For instance, clearly we have $\wt m_1=m_2$, $\wt m_2=m_1$. It  is easy to see that
the fundamental
solution $\{\wt W_k(z)\}$ of the system \eqref{0.1'} is connected with the fundamental solution $\{W_k(z)\}$ of
\eqref{0.1} by the equality
\begin{equation} \label{d2}
\wt W_k(z)=W_k(-z).
\end{equation}
Thus, according to \eqref{0.4} and \eqref{d2} the function
\begin{equation} \label{d3}
\wt \vp(z)=\vp(-z),
\end{equation}
where $\vp$ is the Weyl function of the system \eqref{0.1}, satisfies the inequalities
\begin{align} \label{d4}&
\sum_{k=0}^\infty q(z)^k 
\begin{bmatrix}
 I_{m_2} & \wt \vp(z)^* 
\end{bmatrix}
\wt W_k(z)^* \wt C_k \wt W_k (z)
\begin{bmatrix}
 I_{m_2} \\ \wt \vp(z) 
\end{bmatrix}<\infty \quad (z\in \BC_+).
\end{align}
Therefore, by virtue of  \cite[Definition 3.6]{FKRS14},  the matrix function $\wt \vp(z)$ is the Weyl function (on $\BC_+$)
of the dual system \eqref{0.1'}. Moreover, we see that there is a one to one correspondence
\eqref{d1}, \eqref{d3} between systems \eqref{0.1} and \eqref{0.1'} and their Weyl functions (on $\BC_-$ and $\BC_+$, respectively).
Hence,  \cite[Corollary~4.7]{FKRS14} yields the theorem below.
\begin{Tm}\label{TmUniq} Dirac system \eqref{0.1}, \eqref{0.2} is uniquely recovered from its Weyl function
$\vp(z)$ $(z\in \BC_-)$ introduced by \eqref{0.4}.
\end{Tm}
\paragraph{2.}  In order to consider the case of rational Weyl functions, we introduce generalized B\"acklund-Darboux transformation (GBDT)
of discrete Dirac systems. Each GBDT of the initial discrete Dirac system is determined by a triple $\{A, S_0, \Pi_0\}$ of  parameter matrices.
Here, we take a trivial initial system and
choose $n\in \BN$ $(n>0)$, two $n
\times n$
parameter matrices $A$ ($\det A \not=0$) and
$S_0>0$, and an $n
\times m$ parameter matrix $\Pi_0$ such that
\begin{equation} \label{0.6}
A S_0-S_0 A^*=\I \Pi_0 j \Pi_0^*.
\end{equation}
Define recursively the sequences $\{\Pi_k\}$ and $\{S_k\}$  ($k >0$)
by the
relations
\begin{align} & \label{0.7}
\Pi_{k+1}=\Pi_k+\I A^{-1}\Pi_k j,
\\ &  \label{0.8}
S_{k+1}=S_k+A^{-1}S_k (A^*)^{-1}+A^{-1}\Pi_k
\Pi_k^*(A^*)^{-1}.
\end{align}
From \eqref{0.6}--\eqref{0.8},  the validity of the matrix identity
\begin{equation} \label{0.9}
A S_{r}-S_{r} A^*=\I \Pi_{r} j \Pi_{r}^* \quad
(r \geq 0)
\end{equation}
follows by induction.  
\begin{Dn} \label{adm} The triple $\{A, S_0, \Pi_0\}$, where $\det A \not=0$, $S_0>0$ and 
\eqref{0.6} holds, is called admissible.
\end{Dn}
In view of \eqref{0.8}, for the admissible triple we have $S_k>0$ $(k\geq 0)$. Thus, the
sequence
\begin{equation} \label{0.10}
C_k:=I_m+\Pi_k^*S_k^{-1}\Pi_k-\Pi_{k+1}^*S_{k+1}^{-1}\Pi_{k+1}
\end{equation}
is well-defined. We say that the sequence $\{C_k\}$ is {\it determined} by the
admissible triple $\{A, S_0, \Pi_0\}$.
We will need also the matrix function $w_A$,
which for each $k \geq 0$  is
a so called
transfer matrix function in Lev Sakhnovich form
\cite{SaL1, SaL3, SaSaR} and is defined by the relation
\begin{equation} \label{0.11}
w_A(k,\lambda):=I_m-\I j \Pi_k^*S_k^{-1}(A-\lambda
I_n)^{-1}\Pi_k.
\end{equation}

Now, similar to \cite{KS, FKRS08}, we obtain the  theorem below.
\begin{Tm} \label{FundSol} Let the triple $\{A, S_0, \Pi_0\}$ be admissible
and assume that the recursions \eqref{0.7} and \eqref{0.8} are valid. Then, the
matrices $C_k$ given by \eqref{0.10} $($i.e., determined by $\{A, S_0, \Pi_0\})$
are well-defined and satisfy \eqref{0.2}.
Moreover, in this case the fundamental
solution
$\{W_{k}\}$ of the Dirac system \eqref{0.1}
admits the representation
\begin{equation} \label{0.12}
W_{k}(z)=w_A(k,\,-1/z)\big(I_m +\I z
 j
\big)^{k}w_A(0,\,-1/z)^{-1} \quad (k \geq 0),
\end{equation}
where $w_A$ is defined in \eqref{0.11}.
\end{Tm}
\begin{proof}.
Recall that since $S_0>0$, relation \eqref{0.8} yields by induction that $S_k>0$,
and so the sequence $\{C_k\}$ is well-defined.

Next, formula (\ref{0.12}) easily follows from the 
equality
\begin{equation} \label{0.13}
w_A(k+1, \lambda)\big(I_m - \frac{\I}{\lambda} j
\big)=\big(I_m -
\frac{\I}{\lambda} j C_k \big)w_A(k, \lambda) \quad (k\geq 0),
\end{equation}
which is proved quite similar to the proof of \cite[(2.24)]{FKRS08}  (and so we omit
this proof here).

It remains to prove \eqref{0.2}. The second equality in \eqref{0.2}, that is, $C_kjC_k=j$ follows from
\eqref{0.13} and the equalities
\begin{equation} \label{0.14}
w_A(k,  \lambda)jw_A(k,  \ov \lambda)^*=j,
\end{equation}
which may be found in \cite{SaL1} (see also \cite[(1.84)]{SaSaR}). Indeed, we easily check that
\begin{equation} \label{0.15}
\Big( I_m -\frac{\I}{\lambda} j\Big)j\Big( I_m
+\frac{\I}{\lambda}
j\Big)=\Big( 1 + \frac{1}{\lambda^2} \Big)j, 
\end{equation}
and formulas \eqref{0.13}--\eqref{0.15} imply that
\begin{equation} \label{0.16}
\Big( I_m
-\frac{\I}{\lambda} j C_k \Big)j\Big( I_m + \frac{\I}{\lambda} C_kj\Big)=\Big( 1 +
\frac{1}{\lambda^2}
\Big)j.
\end{equation}
Clearly, the second equality in \eqref{0.2} is immediate from \eqref{0.16}.

Finally, the first equality in \eqref{0.2} is proved in the same way as \cite[Proposition 3.1]{FKRS08}.
\end{proof}
\paragraph{3.}  It is convenient to partition $\Pi_0$ into the $n\times m_i$ blocks $\vt_i$ and to partition $w_A(0,\la)$
in the four blocks of the same orders as for $j$ in \eqref{0.2}:
\begin{align}\label{part}&
\Pi_0=[\vt_1 \quad \vt_2], \quad w_A(0,\lambda)=\left[
\begin{array}{lr}
a(\lambda) & b(\lambda) \\ c(\lambda) & d(\lambda)
\end{array}
\right].
\end{align}
\begin{Tm} \label{TmDirect}
Let a sequence
$\{C_k\}$ and so Dirac system \eqref{0.1}, \eqref{0.2} be determined
by some admissible triple $\{A,S_0, \Pi_0\}$.
Then, the unique Weyl function of this system
is given by the formula
\begin{equation} \label{1.1}
\vp(z)=-\I z\vt_1^*S_0^{-1}(I_n+zA^{\times})^{-1}\vt_2,
\quad A^{\times}=A+\I \vt_2 \vt_2^*S_0^{-1}.
\end{equation}
\end{Tm}
\begin{proof}. 
Recall the definition \eqref{0.4}  of the Weyl function $\vp(z)$, where $q(z)=(1+|z|^2)^{-1}$.
First, let us show that the summation formula
\begin{equation} \label{1.2}
\sum_{k=0}^r
q(z)^kW_k(z)^*C_kW_k(z)=\frac{\I (1+|z|^2)}{(
 \ov{z}-z
)}\Big(q(z)^{r+1}W_{r+1}(z)^*jW_{r+1}(z)-j\Big)
\end{equation}
is valid.
Indeed, according to (\ref{0.1}) and (\ref{0.2}) we
have
\begin{align}\nn
W_{k+1}(z)^*jW_{k+1}(z)&=W_{k}(z)^*\big(
I_m -\I \ov{z}
 C_k j \big)j\big( I_m
+\I z j
C_k \big)W_{k}(z)
\\ & \nn
=q(z)^{-1}W_k(z)^*jW_k(z)+{\I(
z -\ov
z)}W_k(z)^*C_kW_k(z),
\end{align}
that is,
\begin{align}\label{1.3}
q(z)^kW_k(z)^*C_kW_k(z)=&\frac{\I q(z)^{k-1}}{(\ov{z}-z)}
\\ \nn & \times
\left(q(z)W_{k+1}(z)^*jW_{k+1}(z)-W_k(z)^*jW_k(z)\right),
\end{align}
and \eqref{1.2} is immediate from \eqref{1.3}. 

Next, we will need the inequality
\begin{align}\label{1.4}
w_A\left(k, -\frac{1}{{z}}\right)^*j w_A\left(k, -\frac{1}{z}\right) \leq j \quad (z \in \BC_-),
\end{align}
which together with \eqref{0.14} follows from a more general formula (see, e.g., \cite[(1.88)]{SaSaR}) of the form
\begin{align}\label{1.5}
w_A\left(k, {\la}\right)^*j w_A\left(k, \la\right) =j -\I (\la - \ov{\la})\Pi_k^*(A^*-\ov{\la}I_n)^{-1} S_k^{-1} (A-\la I_n)^{-1}\Pi_k.
\end{align}
Formulas  \eqref{0.12} and \eqref{1.4} yield (in $\BC_-$) the inequality
\begin{align}\nn
W_{r+1}(z)^*jW_{r+1}(z) \leq &\big(w_A(0,\,-1/z)^{-1} \big)^* \big(I_m-\I \ov{z}j)^{r+1}j  
\\ \label{1.6} & \times
 \big(I_m +\I z
 j\big)^{r+1}w_A(0,\,-1/z)^{-1}.
\end{align}
Setting
\begin{align}\label{1.7}
\vp(z)=b(-1/z)d(-1/z)^{-1}
\end{align}
and taking into account \eqref{part} and \eqref{1.7}, we derive
\begin{align}
 \nn
 \big(I_m +\I z
 j\big)^{r+1}w_A(0,\,-1/z)^{-1}\begin{bmatrix}\vp(z) \\ I_{m_2} \end{bmatrix}&=
 \big(I_m +\I z
 j\big)^{r+1}\begin{bmatrix}0 \\ I_{m_2} \end{bmatrix}d(-1/z)^{-1}
 \\ \label{1.8}  &
 =(1-\I z)^{r+1}\begin{bmatrix}0 \\ I_{m_2} \end{bmatrix}d(-1/z)^{-1}.
\end{align}
It is immediate from \eqref{1.6} and \eqref{1.8} that
\begin{align} \label{1.9}
\begin{bmatrix}\vp(z)^* & I_{m_2} \end{bmatrix}W_{r+1}(z)^*jW_{r+1}(z) \begin{bmatrix}\vp(z) \\ I_{m_2} \end{bmatrix}
\leq 0 \quad (z \in \BC_-).
\end{align}
For $\vp(z)$ given by \eqref{1.7}, relations  \eqref{1.2} and \eqref{1.9} imply that \eqref{0.4} holds, and so this $\vp(z)$
is the Weyl function. (We did not discuss the singularities of $d(-1/z)$ and $d(-1/z)^{-1}$ but $\vp(z)$ is analytic in $\BC_-$ because
it is meromorphic and it is the Weyl function.)

It remains to show that the right-hand sides of \eqref{1.1} and \eqref{1.7} coincide. By virtue of
\eqref{0.11} and \eqref{part}, using inversion
formula from system theory (see, e.g., \cite[Appendix B]{SaSaR} and references therein),
 we obtain
\begin{align}\nn
b(\la)d(\la)^{-1}&=-\I \vt_1^*S_0^{-1}(A-\la I_n)^{-1}\vt_2\big(I_{m_2}+\I \vt_2^* S_0^{-1}(A-\la I_n)^{-1} \vt_2\big)^{-1}
\\ \nn &
=-\I \vt_1^*S_0^{-1}(A-\la I_n)^{-1}\vt_2\big(I_{m_2}-\I \vt_2^* S_0^{-1}(A^{\times}-\la I_n)^{-1} \vt_2\big),
\end{align}
where $A^{\times}  =A+\I \vt_2\vt_2^* S_0^{-1}$. Since $\I \vt_2\vt_2^* S_0^{-1}=A^{\times} -A=(A^{\times}- \la I_n) -(A-\la I_n)$,
we essentially simplify the right-hand side in the formula above:
\begin{align}\label{1.10}&
b(\la)d(\la)^{-1}=-\I\vt_1^* S_0^{-1}(A^{\times}-\la I_n)^{-1} \vt_2.
\end{align}
Hence, the right-hand sides of \eqref{1.1} and \eqref{1.7}, indeed, coincide.
\end{proof}

\paragraph{4.} We note that the Weyl function $\vp(z)$ in \eqref{1.1} is rational and contractive on $\BC_-$.
Moreover, $\vp(-1/z)$ is strictly proper rational and contractive.
It is well-known (see, e.g., \cite{KFA, LR}) that each strictly proper rational $m_1\times m_2$ matrix function
$\psi(z)$ admits a representation (so called {\it realization})
\begin{equation}
\label{1.11} \psi(z)=\clc(z I_n- \cla)^{-1}\clb,
\end{equation}
where $\cla$ is an $n \times n$ matrix, $\clc$ is an $m_1 \times n$ matrix and $\clb$ is an $n\times m_2$ matrix.
Further in the text we assume that the realization  \eqref{1.11} is a {\it minimal realization}, that is, the value of $n$
in \eqref{1.11} is minimal (among the corresponding values in different realizations of $\psi$).
The following proposition is immediate from  \cite[Lemma 3.1]{ALS16}
(and is based on several theorems from \cite{LR}, see the details in \cite{ALS16}). 
\begin{Pn}\label{PnRic} Assume that a strictly proper rational $m_1\times m_2$ matrix function $\psi(z)$ is contractive on $\BC_-$
and that \eqref{1.11} is its minimal realization. 
Then, there is a unique Hermitian
solution $X$ of the Riccati equation 
\begin{align} \label{1.12}&
 X\clb \clb^*X-\I(\cla^*X-X\cla )+\clc^*\clc =0.
\end{align}
such that 
the relation
\begin{align} &       \label{1.13}
\s(\cla-\I \clb\clb^* X)\subset (\BC_+\cup \BR)
\end{align} 
holds. Moreover, this solution $X$  is positive.
\end{Pn}
Next, we give an explicit procedure of solving the inverse problem to recover  Dirac system from its Weyl function.
\begin{Tm}\label{TmIP} Let $\vp(z)$ be a rational $m_1 \times m_2$ matrix function such that $\psi(z)=\vp(-1/z)$ 
is a strictly proper rational matrix function, which  is contractive on $\BR$ and has no poles on $\BC_-$.  Assume that \eqref{1.11} is a minimal realization
of $\psi$ and that $X>0$ is a solution of \eqref{1.12}. 

Then, $\vp(z)$ is the Weyl function of the Dirac system \eqref{0.1}, \eqref{0.2}, the potential $\{C_k\}$ of which
is determined by the admissible triple
\begin{align} 
\label{1.14}&
A=\cla -\I \clb \clb^*X,  \quad S_0=X^{-1}, \quad \vt_1=\I X^{-1}\clc^*, \quad \vt_2=\clb.
\end{align}
\end{Tm}
\begin{proof}.
Since $\psi(z)$ is contractive on $\BR$ and has no poles on $\BC_-$, it is contractive
on $\BC_-$. Thus, according to Proposition \ref{PnRic} a positive definite solution $X$ of \eqref{1.12} exists.
In view of \eqref{1.14}, choosing $X>0$ we have $S_0>0$. Moreover, relations \eqref{1.12} and \eqref{1.13}
yield the equality
\begin{align} \label{1.15}&
 \vt_2\vt_2^*+\I\big((A+\I\vt_2\vt_2^*S_0^{-1})S_0-S_0(A+\I\vt_2\vt_2^*S_0^{-1})^*\big)+\vt_1\vt_1^* =0,
\end{align}
which is equivalent to \eqref{0.6}.  Hence, the triple $\{A,S_0, \Pi_0\}$ is admissible.

It remains to show that for the Weyl function $\vp(z)$ of the Dirac system (determined by this triple),
the function $\psi(z)=\vp(-1/z)$ coincides with  $\psi(z)$ admitting the realization \eqref{1.11}.
Taking into account Theorem \ref{TmDirect} and equalities \eqref{1.14}, we see that $\psi(z)$ determined by our triple has the form
\begin{align} \label{1.16}&
\psi(z)=\I \vt_1^*S_0^{-1}(zI_n-\cla)^{-1}\vt_2=\clc(zI_n-\cla)^{-1}\clb ,
\end{align}
and the right-hand sides of \eqref{1.11} and \eqref{1.16}, indeed, coincide.
\end{proof} 
\section{Verblunsky-type coefficients and asymptotics of the potentials}   \label{Verb}
\setcounter{equation}{0}
Recall that the matrices $C_k$ from the potential (sequence) $\{C_k\}$ are positive definite and $j$-unitary
(i.e., they satisfy \eqref{0.2}). According to \cite[Proposition~2.4]{FKRS14} it means that they admit representations
\begin{align}\label{v1}&
C_k= {\cal D}_k  H_k, \quad   {\mathcal D}_k:= {\mathrm{diag}}\Big\{
\big(
I_{m_1}-  \rho_k  \rho_k^* \big)^{-\frac{1}{2}}, \, \,
\big(I_{m_2}- \rho_k^* \rho_k\big)^{-\frac{1}{2}}\Big\},
\\ \label{v2}&
H_k:= \left[
\begin{array}{cc}
I_{m_1} &  \rho_k \\   \rho_k^* & I_{m_2} 
\end{array}
\right] \quad (\rho_k^*\rho_k <I_{m_2}).
\end{align}
Here, the $m_1 \times m_2$ matrices $\rho_k$ are so called Verblunsky-type coefficients,
which were studied in detail in \cite{ALSverb}. It is well-known (see, e.g., \cite{DFK}) that ${\cal D}_k  H_k=H_k\cld_k$.
Clearly, $\rho_k^*\rho_k <I_{m_2}$ yields $\rho_k\rho_k^* <I_{m_1}$ and vice versa.

In this section, we show that 
\begin{align} \label{v3}&
\lim_{k\to \infty} \begin{bmatrix} I_{m_1} & 0\end{bmatrix}C_k\begin{bmatrix} I_{m_1} \\ 0\end{bmatrix}=I_{m_1},
\end{align}
and so $\rho_k \to 0$ and $C_k\to I_m$. More precisely, we prove the following statement.
\begin{Tm}\label{TmAs} Let the triple $\{A, S_0, \Pi_0\}$ be admissible and assume that $-\I \not\in\s(A)$.
Then, for the potential
 $\{C_k\}$ $($of the Dirac system \eqref{0.1}$)$ determined by this triple the asymptotic relations 
\begin{align} \label{v4}&
\lim_{k\to \infty} \rho_k=0, \quad \lim_{k\to \infty} C_k=I_m
\end{align} 
 are valid.
 \end{Tm}
\begin{proof}.  Consider the equality
\begin{align} \nn &
S_{k+1}-(I_n+\I A^{-1})S_k\big(I_n-\I (A^*)^{-1}\big)
\\ \label{v5} 
&=S_{k+1}-S_k- A^{-1}S_k (A^*)^{-1}+\I A^{-1}(AS_k-S_kA^*)(A^*)^{-1}.
\end{align} 
Using \eqref{0.8} and \eqref{0.9}, we rewrite \eqref{v5}:
\begin{align}  \label{v6}  &
S_{k+1}-(I_n+\I A^{-1})S_k\big(I_n-\I (A^*)^{-1}\big)=A^{-1}\Pi_k(I_m-j)\Pi_k^*(A^*)^{-1}.
\end{align} 
Now, we partition $\Pi_k$ and, taking into account \eqref{0.7} and \eqref{part}, write it down in the form
\begin{align}  \label{v7}  &
\Pi_k=\begin{bmatrix}(I_n+\I A^{-1})^k\vt_1 & (I_n-\I A^{-1})^k \vt_2 \end{bmatrix}.
\end{align} 
In view of \eqref{v6} and \eqref{v7}, setting
\begin{align}  \label{v8}  &
R_{r}:=(I_n+\I A^{-1})^{-r}S_r\big(I_n-\I (A^*)^{-1}\big)^{-r}
\end{align} 
we have
\begin{align} \nn
R_{k+1}-R_k=& 2(I_n+\I A^{-1})^{-k-1}A^{-1}(I_n-\I A^{-1})^k \vt_2\vt_2^* \big((I_n-\I A^{-1})^k\big)^*\big(A^{-1}\big)^*
\\  \label{v9}  &
\times
\big((I_n+\I A^{-1})^{-k-1}\big)^*\geq 0.
\end{align} 
Since $R_0=S_0>0$, relations \eqref{v9} imply that there is a
limit
\begin{align}  \label{v10}  &
\lim_{k\to \infty}R_k^{-1}=\vk_R\geq 0.
\end{align} 
On the other hand, from \eqref{v7} and \eqref{v8} we derive
\begin{align} \label{v11}&
 \begin{bmatrix} I_{m_1} & 0\end{bmatrix}\Pi_k^*S_k^{-1}\Pi_k\begin{bmatrix} I_{m_1} \\ 0\end{bmatrix}=\vt_1^*R_k^{-1}\vt_1,
\end{align}
and so \eqref{v10} yields
\begin{align} \label{v12}&
\lim_{k \to \infty} \begin{bmatrix} I_{m_1} & 0\end{bmatrix}\Pi_k^*S_k^{-1}\Pi_k\begin{bmatrix} I_{m_1} \\ 0\end{bmatrix}=\vt_1^* \vk_R\vt_1,
\end{align}
The definition \eqref{0.10} of $C_k$ and the existence of the limit in \eqref{v12} show that \eqref{v3} holds. It is easy to see that the first
equality in \eqref{v4} follows from \eqref{v1}--\eqref{v3}.  Finally, the second equality in \eqref{v4} is immediate from \eqref{v1}, \eqref{v2}
and the first
equality in \eqref{v4}.
\end{proof}
\begin{Rk} According to Theorems \ref{TmDirect}, \ref{TmIP} and \ref{TmUniq}, and to Proposition \ref{PnRic},
given a potential $\{C_k\}$ determined by some admissible triple we may recover another admissible triple $\{A, S_0, \Pi_0\}$, which
determines the same sequence $\{C_k\}$ and has additional property $\s(A)\subset (\BC_+\cup \BR)$.
Namely, we construct first  the Weyl function using the initial triple and the procedure from Theorem \ref{TmDirect}. Next,
we recover another admissible triple $\{A, S_0, \Pi_0\}$ such that $\s(A)\subset (\BC_+\cup \BR)$
in the process of solving inverse problem.

Thus, we may assume $\s(A)\subset (\BC_+\cup \BR)$ without loss of generality, and so the condition
$-\I \not\in \s(A)$ in Theorem \ref{TmAs} may be deleted.
\end{Rk}
We note that in the case of $\{C_k\}$ determined by some admissible triple, Verblunsky-type coefficients
may be expressed explicitly. Indeed, in view of \eqref{v1} and \eqref{v2} we have
\begin{align} \label{v13}&
\rho_k=\left( \begin{bmatrix} I_{m_1} & 0\end{bmatrix}C_k\begin{bmatrix} I_{m_1} \\ 0\end{bmatrix}\right)^{-1}
 \begin{bmatrix} I_{m_1} & 0\end{bmatrix}C_k\begin{bmatrix}0 \\ I_{m_2} \end{bmatrix}.
\end{align}
Hence, taking into account \eqref{0.10} and \eqref{v11} we derive
\begin{align}\nn
\rho_k=&\left(I_{m_1}+\vt_1^*R_k^{-1}\vt_1-\vt_1^*R_{k+1}^{-1}\vt_1\right)^{-1}
\\  \label{v14}& \times
 \begin{bmatrix} I_{m_1} & 0\end{bmatrix}(\Pi_k^*S_k^{-1}\Pi_k-\Pi_{k+1}^*S_{k+1}^{-1}\Pi_{k+1})\begin{bmatrix}0 \\ I_{m_2} \end{bmatrix}.
\end{align}
\section{Stability of solving inverse problem}\label{stab}
\setcounter{equation}{0}
It is easy to see that the procedure (given in Theorem \ref{TmIP}) to recover system \eqref{0.1}, \eqref{0.2}  
consists from two steps. The first step is the construction of $X>0$ and the second step is the
construction of the potential $\{C_k\}$ using this $X$. 

We start with  the matrix function $\vp(z)$ such that $\psi(z)=\vp(-1/z)$ is a strictly proper rational $m_1 \times m_2$ matrix function,
which is contractive  on $\BC_-$. More precisely, we start with a minimal realization \eqref{1.11} of $\psi$ (or, equivalently,
with the triple $\{ \cla,  \clb,  \clc\}$) and consider the stability in recovery of $X>0$ satisfying additional condition \eqref{1.13}. The existence and uniqueness
of  $X>0$ satisfying \eqref{1.13} follows from Proposition \ref{PnRic}.

\begin{Dn}\label{stabX}
By $\clg_n$ we denote the class of triples $\{\wt \cla, \wt \clb, \wt \clc\}$ which determine minimal realizations $\wt \psi(z)=\wt \clc(zI_n-\wt \cla)^{-1}\wt \clb$
of $m_1\times m_2$ matrix functions $\wt \psi(z)$ contractive on $ \BC_-$. 
 
The recovery of $X>0$ satisfying \eqref{1.12}, \eqref{1.13}  from the minimal realization \eqref{1.11}  of $\psi(z)$ $($where $\{ \cla,  \clb,  \clc\}\in \clg_n)$ is called  stable
if for any $\ve>0$ there is $\delta >0$ such that for each $\{\wt \cla, \wt \clb, \wt \clc\}$, satisfying conditions
\begin{align} &       \label{1.17}
\{\wt \cla, \wt \clb, \wt \clc\}\in \clg_n, \quad \|\cla -\wt \cla\|+ \| \clb -\wt \clb\|+\|\clc- \wt \clc\|<\delta,
\end{align} 
there is a solution $\wt X=\wt X^*$ of the equation 
\begin{align} &       \label{1.18}
 \wt X \wt \clb \wt \clb^* \wt X-\I(\wt \cla^* \wt X-\wt X\wt \cla )+\wt \clc^* \wt \clc =0
\end{align} 
in the neighbourhood $\| X-\wt X\|<\ve$ of $X$.
\end{Dn}
The stability of the recovery of $X$ follows (similar to the case of continuous Dirac system) from \cite[Theorem 3.3]{ALS16} based on
 \cite[Theorem 4.4]{RaRo}. Namely, applying  \cite[Theorem 3.3]{ALS16} to the triples $\{-\cla,\clb, -\clc\}$ and $\{-\wt \cla,\wt \clb, -\wt \clc\}$
we obtain our next statement.
\begin{Pn}\label{PnStabD} The recovery of $X>0$, satisfying \eqref{1.12}, \eqref{1.13}  from the minimal realization \eqref{1.11} $($with $\{ \cla,  \clb,  \clc\}\in \clg_n)$
  is stable.
\end{Pn}
 \begin{Rk}\label{Rkclg} Note that $($according to \cite[Theorem 4.4]{RaRo}$)$ we may consider  a wider than $\clg_n$ class of  perturbed triples 
$\{\wt \cla, \wt \clb, \wt \clc\}$, that is, such perturbed triples that \eqref{1.18} has a Hermitian solution $\wt X=\wt X^*$.
 \end{Rk}
 Recall that given the triple $\{\cla, \clb, \clc\}$ and $X>0$ we construct the matrices $A, \, S_k, \, R_k, \ldots$
 For the matrices constructed in a similar way in the case of the triple $\{\wt \cla, \wt \clb, \wt \clc\}$ and of $\wt X>0$
 satisfying
\begin{align} &       \label{Ric'}
\wt X \wt \clb \wt \clb^* \wt X-\I(\wt \cla^*\wt X-\wt X\wt \cla )+\wt \clc^*\wt\clc =0,
\end{align}   
we use the notations with ``tilde":  $\wt A, \, \wt S_k, \, \wt R_k, \ldots$

The stability of the second step of solving inverse problem one can prove under additional condition $\vk_R=0$ or, equivalently,
\begin{align} &       \label{1.19}
\lim_{k\to \infty}R_k=+\infty,
\end{align} 
which means that all the eigenvalues of $R_k$ tend to infinity. Unlike the skew-self-adjoint case \cite{FKRS17}, the equality
\eqref{1.19} is not fulfilled automatically. 

Sufficient condition of stability may be expressed also in terms of matrices $Q_r$, which are introduced
by the relations
\begin{align}  \label{v8'}  &
Q_{r}:=(I_n-\I A^{-1})^{-r}S_r\big(I_n+\I (A^*)^{-1}\big)^{-r} .
\end{align} 
Clearly, we assume in \eqref{v8'} that $\I \not\in \s(A)$.
Similar to  the equality \eqref{v6}, from \eqref{0.8} and \eqref{0.9} we have
\begin{align}  \label{v6'}  &
S_{k+1}-(I_n-\I A^{-1})S_k\big(I_n+\I (A^*)^{-1}\big)=A^{-1}\Pi_k(I_m+j)\Pi_k^*(A^*)^{-1}.
\end{align} 
Hence, taking into account \eqref{v7}  (in analogy with the relation \eqref{v9} for $R_r$) we derive
\begin{align} \nn
Q_{k+1}-Q_k=& 2(I_n-\I A^{-1})^{-k-1}A^{-1}(I_n+\I A^{-1})^k \vt_1\vt_1^* \big((I_n+\I A^{-1})^k\big)^*\big(A^{-1}\big)^*
\\  \label{v9'}  &
\times
\big((I_n-\I A^{-1})^{-k-1}\big)^*\geq 0.
\end{align} 
Since $Q_0=S_0>0$, relations \eqref{v9'} imply that there is a
limit
\begin{align}  \label{v10'}  &
\lim_{k\to \infty}Q_k^{-1}=\vk_Q\geq 0.
\end{align} 
Moreover, \eqref{v7} and \eqref{v8'} yield
\begin{align} \label{v12'}&
\lim_{k \to \infty} \begin{bmatrix}0 & I_{m_1} \end{bmatrix}\Pi_k^*S_k^{-1}\Pi_k\begin{bmatrix}0 &  I_{m_1} \end{bmatrix}=\vt_2^* \vk_Q\vt_2.
\end{align}
Formula \eqref{v12'} implies that
\begin{align} \label{v3'}&
\lim_{k\to \infty} \begin{bmatrix}0 & I_{m_2} \end{bmatrix}C_k\begin{bmatrix}0 & I_{m_2} \end{bmatrix}=I_{m_2},
\end{align}
which gives another way to prove Theorem \ref{TmAs}. The cases when \eqref{1.19} or the equality
\begin{align} &       \label{1.19'}
\lim_{k\to \infty}Q_k=+\infty
\end{align} 
hold are considered in the stability theorem below. (Recall that the sequence $\{R_k\}$ is given by \eqref{v8} or, equivalently, by \eqref{v9} 
together with \eqref{1.14} and $R_0=S_0.)$
In Proposition \ref{PnEx} at the end of this section we present a wide
class, where  \eqref{1.19'} is valid. 
\begin{Tm} \label{TmdsaStab} Consider the procedure $($from Theorem \ref{TmIP}$)$ 
of the unique recovery of the  potential $\{C_k\}$ of  the discrete self-adjoint Dirac system \eqref{0.1}, \eqref{0.2}
 from a minimal
realization \eqref{1.11}, where $\psi(z)=\vp(-1/z)$ and $\vp(z)$ is the Weyl function of the system \eqref{0.1}, \eqref{0.2}.
Assume that $X$ in this procedure is chosen so that \eqref{1.13} holds $($which is always possible$)$.
Assume also that either the sequence $\{R_k\}$
satisfies \eqref{1.19} or $\I \not\in \s(A)$ and the sequence $\{Q_k\}$
satisfies \eqref{1.19'}.

Then, this procedure of the recovery of the  potential $\{C_k\}$ is stable in the class of the triples from $ \clg_n$.
\end{Tm}
\begin{proof}. The recovery of $X>0$ satisfying \eqref{1.12}, \eqref{1.13} is possible according to Proposition \ref{PnRic}
and is stable according to Proposition \ref{PnStabD}. 

Now, in order to show that the recovery of $\{C_k\}$ is stable  under condition \eqref{1.19}, we choose
some small $\wh \ve >0$ and such a large $N>0$ and a small  neighbourhood of $\{\cla, \clb, \clc\}$ that $\|R_k^{-1}\|<\wh \ve$
and $\|\wt R_k^{-1}\|<2 \wh \ve$ for $X>0$ satisfying \eqref{1.12}, \eqref{1.13}, for $k>N$, and for the matrices
$\wt X>0$
satisfying \eqref{Ric'}
(where the triples $\{\wt \cla, \wt \clb, \wt \clc\} \in \clg_n$ belong to the mentioned above neighbourhood of $\{\cla, \clb, \clc\}$ 
and $\wt X$ are those solutions of \eqref{Ric'} which belong to the  neighbourhood of  $X$).
Here, we use the fact that the sequence $\{\wt R_k\}$ is monotonically increasing and if $\wt R_{r_0}$ is sufficiently large,
then $\wt R_r$ ($r>r_0$) is sufficiently large as well. 

In view of  \eqref{0.10} and \eqref{v11}, we see that for sufficiently small
$\wh \ve$ the matrices
\begin{align} \label{1.20}&
\begin{bmatrix} I_{m_1} & 0\end{bmatrix}C_k\begin{bmatrix} I_{m_1} \\ 0\end{bmatrix}, \quad \begin{bmatrix} I_{m_1} & 0\end{bmatrix}
\wt C_k\begin{bmatrix} I_{m_1} \\ 0\end{bmatrix},
\end{align}
 are sufficiently close to $I_{m_1}$. This, in turn, means that (in view of \eqref{v1} and \eqref{v2})
the matrices $\rho_k$, $\wt \rho_k$ are sufficiently small, and so $C_k$ and $\wt C_k$ are sufficiently
close to $I_m$. Therefore, for any $\ve>0$ we may choose $\wh \ve$ such that 
$$\|C_k-\wt C_k\|<\ve
\quad {\mathrm{for \,\, all}} \quad k>N(\wh \ve).$$

Moreover, for any $\ve >0$  we may choose a  neighbourhood of $X$ and of $\{\cla, \clb, \clc\}$
such that for $\{\wt \cla, \wt \clb, \wt \clc\}$ from this  neighbourhood
the inequalities
$$
\|C_k-\wt C_k\|<\ve \quad (0 \leq k \leq N(\wh \ve))
$$
are valid as well. Thus, the recovery of $\{C_k\}$ is stable, indeed.

The stability of the  recovery of $\{C_k\}$ under condition \eqref{1.19'} is proved in a similar way.
\end{proof}
Now, consider the case when $A$ is similar to a diagonal matrix $D$ ($A$ is diagonalisable):
\begin{align} \label{1.22}&
A=UDU^{-1}.
\end{align}
Relations \eqref{1.13}, \eqref{1.14} and \eqref{1.22} yield $\s(D)\in (\BC_+\cup \BR)$ or, equivalently:
\begin{align} \label{1.23}&
\I(D^*-D)\geq 0.
\end{align}
\begin{Pn}\label{PnEx} Let the sequence $\{Q_k\}$ be given by \eqref{v8'}, where $A$ and $\{S_k\}$
are constructed using the procedure from Theorem \ref{TmdsaStab}, $A$ is diagonalisable
$($i.e., the representation \eqref{1.22} holds$)$ and $\I \not\in \s(A)$.
Then, \eqref{1.19'} is valid.
\end{Pn}
\begin{proof}. According to \eqref{v9'} we have
 \begin{align} \label{1.30} &
 Q_{k+n}-  Q_k
 =  2(A-\I I_n)^{-n-k}(A+\I I_n)^{k}F
 (A^*-\I I_n)^{k}
 (A^*+\I I_n)^{-n-k},
\\ & \label{1.31}
F:=\sum_{\ell=1}^n (A-\I I_n)^{n-\ell}(A+\I I_n)^{\ell-1} \vt_1  \vt_1^*(A^*-\I I_n)^{\ell - 1}(A^*+\I I_n)^{n-\ell},
 \end{align} 
 where $F$ does not depend on $k$. Let us show that $F$ is strictly positive, that is, $F>0$.
Indeed, it is easy to see (more details are given in the similar part of the proof of \cite[Proposition 4.10]{FKRS17}) that
$$\spa\bigcup_{\ell=1}^n(A-\I I_n)^{n-\ell}(A+\I I_n)^{\ell-1} \vt_1 = \spa\bigcup _{\ell=1}^n A^{\ell-1} \vt_1,
$$ 
and so we need only to prove that the pair $\{A, \vt_1\}$ is controllable.

Since the realization \eqref{1.11} is minimal, the pair $\{\cla^*, \clc^*\}$ is controllable.
In view of \eqref{1.14}, the controllability of  the pair $\{X^{-1}\cla^*X,  \vt_1\}$ follows from the
controllability of $\{\cla^*, \clc^*\}$. Hence, the equality
\begin{align} \label{1.32}&
X^{-1}\cla^*X=A-\I \vt_1\vt_1^* X
\end{align}
(which we derive below) implies that the pair $\{A, \vt_1\}$ is controllable as well.

Finally, using \eqref{1.14} we rewrite \eqref{0.6} in the form
$$AX^{-1}-X^{-1}A^*=\I (\vt_1\vt_1^*-\vt_2\vt_2^*).$$
This yields in turn that
$X^{-1}A^*X=A+\I \clb\clb^*X-\I\vt_1\vt_1^*X$. Applying now the first equality
in \eqref{1.14}, we obtain \eqref{1.32}, and so $\{A, \vt_1\}$ is controllable and
the inequality $F>0$ is proved.

Next, we show that
\begin{align} \label{1.33}&
(D-\I I_n)^{-1}(D+\I I_n)\big((D-\I I_n)^{-1}(D+\I I_n)\big)^* \geq I_n.
\end{align}
The inequality \eqref{1.33} is equivalent to the inequality
$$(D+\I I_n)(D^*-\I I_n) \geq (D-\I I_n)(D^*+\I I_n),$$
which follows from \eqref{1.23}.

Now, formula \eqref{1.30}, representation \eqref{1.22} and inequalities $F>0$ and \eqref{1.33}
imply that
\begin{align} \label{1.34}&
Q_{k+n}-  Q_k\geq \ve I_n
\end{align}
for some $\ve >0$, which does not depend on $k$. The asymptotics \eqref{1.19'} is immediate
from \eqref{1.34}.
\end{proof}

\bigskip

\noindent{\bf Acknowledgments.}
 {This research   was supported by the
Austrian Science Fund (FWF) under Grant  No. P29177.}

\section{Appendix}
\begin{proof} of Theorem \ref{TmWf}. 
It is easy to see that
\begin{align}\label{ap0}&
(I_m+\I \ov{z} jC_k )^*j(I_m+\I z jC_k )=(1+z^2)j,
\end{align}
and so both $(I_m+\I z jC_k )$ and $W_r(z)=\prod_{k=0}^{r-1}(I_m+\I z jC_k )$ are invertible for $z\not=\pm\I$.
Now, let us consider the sets $\cln_r$ of the  linear fractional transformations
\begin{align}\label{ap1}&
\vp_r(z, \clp)=\begin{bmatrix}
I_{m_1} & 0
\end{bmatrix}W_{r}(z)^{-1}\clp(z)\Big(\begin{bmatrix}
0 & I_{m_2}
\end{bmatrix}W_{r}(z)^{-1}\clp(z)\Big)^{-1},
\end{align}
where $\clp(z)$ are nonsingular $m \times m_2$ matrix functions with property-$j$. 
That is, $\clp(z)$ are meromorphic on $\BC_-$ matrix functions such that the inequalities
\begin{align}\label{ap2}&
\clp(z)^*\clp(z)>0, \quad \clp(z)^*j\clp(z) \leq 0
\end{align}
hold for all the points in $\BC_-$ (excluding, possibly, discrete sets of points). 
The sets $\cln_r$  are well-defined because the inequality
\begin{align}\label{ap1+}&
\det\Big(\begin{bmatrix}
0 & I_{m_2}
\end{bmatrix}W_{r}(z)^{-1}\clp(z)\Big)\not=0
\end{align}
follows from \eqref{ap2}.
Indeed, since relations \eqref{0.2} and
\eqref{Prop2.2} yield
\begin{align}\label{ap3}&
(I_m+\I z jC_k )^*j(I_m+\I z jC_k )=(1+|z|^2)j+\I(z-\ov{z})C_k\geq \wt q(z) j, 
\\ \label{ap4}&
\wt q(z):=1+|z|^2+\I (z-\ov{z})>0,
\end{align}
we have 
\begin{align}\label{ap5}&
W_r(z)^*j W_r(z)\geq \wt q(z)^r j, \quad  {\mathrm{i.e.,}} \quad \big(W_r(z)^{-1}\big)^*j W_r(z)^{-1}\leq \wt q(z)^{-r} j.
\end{align}
Thus, the inequalities
\begin{align}\label{ap6}&
\clp(z)^* \big(W_r(z)^{-1}\big)^*j W_r(z)^{-1}\clp(z)\leq 0, \quad \begin{bmatrix}
0 & I_{m_2}
\end{bmatrix}j \begin{bmatrix}
0 \\ I_{m_2}
\end{bmatrix}<0
\end{align}
are valid, and \eqref{ap1+} is immediate from \cite[Proposition 1.43]{SaSaR}.

In view of \eqref{ap1} we have
\begin{align}\label{ap7}&
\vp_{r+1}(z, \clp)=\begin{bmatrix}
I_{m_1} & 0
\end{bmatrix}W_{r}(z)^{-1}\wt \clp(z)\Big(\begin{bmatrix}
0 & I_{m_2}
\end{bmatrix}W_{r}(z)^{-1}\wt \clp(z)\Big)^{-1},
\end{align}
where
\begin{align}\label{ap8}&
\wt \clp(z)=(I_m+\I z jC_r )^{-1}\clp(z).
\end{align}
Relations \eqref{ap3}, \eqref{ap4} and \eqref{ap8} imply that
\begin{align}\label{ap9}&
\wt \clp(z)^*j\wt \clp(z)\leq 0.
\end{align}
Compare \eqref{ap1}, \eqref{ap2} with \eqref{ap7}, \eqref{ap9} to see that the sets (Weyl disks) $\cln_r$
are embedded:
\begin{align}\label{ap10}&
\cln_{r+1} \subseteq \cln_r.
\end{align}
Clearly, formulas \eqref{ap7}--\eqref{ap9} remain valid when we put there $r=0$.
For that case, we partition $\wt \clp$ and (in view of \eqref{0.3}) we rewrite \eqref{ap7} in the form
\begin{align}\label{ap11}&
\vp_1(z,\clp)=\wt \clp_1(z)\wt \clp_2(z)^{-1}, \quad \wt \clp=: \begin{bmatrix}
\wt \clp_1 \\ \wt \clp_2
\end{bmatrix},
\end{align}
where (according to \eqref{ap1+} with $r=1$) we have $\det \wt \clp_2(z)\not=0$.
It follows from \eqref{ap9} and \eqref{ap11} that the functions from $\cln_1$ are contractive. 
Hence, \eqref{ap10} implies that all the functions $\vp_r(z, \clp)$ 
given by \eqref{ap1} are analytic and contractive in $\BC_-$.

Next, using Montel's theorem and arguments from the Step 1 in the proof of \cite[Theorem 3.8]{FKRS14}
one may easily show that there is an analytic and contractive in $\BC_-$ matrix function $\vp_{\infty}(z)$
such that
\begin{align}\label{ap12}&
\vp_{\infty} \in \bigcap_{r \geq 1} \cln_r .
\end{align}
(We note the functions $ \begin{bmatrix}
I_{m_1} \\ \vp
\end{bmatrix}$ in the proof of \cite[Theorem 3.8]{FKRS14}
should be substituted by $ \begin{bmatrix}
\vp \\ I_{m_2} 
\end{bmatrix}$ 
for our case of Weyl functions in $\BC_-$.) Taking into account \eqref{ap1} and \eqref{ap12}
we write the representations
\begin{align}\label{ap13}&
 \begin{bmatrix}
\vp_{\infty}(z) \\ I_{m_2} 
\end{bmatrix}=W_{r+1}(z)\clp(z,{r+1}) \quad (r \geq 0),
\end{align}
where $\clp(z,{r+1})$ are nonsingular with property-$j$. Using the summation formula \eqref{1.2} and representation
\eqref{ap13}, we derive
\begin{align}\label{ap14}&
 \begin{bmatrix}
\vp_{\infty}^* & I_{m_2} 
\end{bmatrix}\sum_{k=0}^r
q(z)^kW_k(z)^*C_kW_k(z) \begin{bmatrix}
\vp_{\infty} \\ I_{m_2} 
\end{bmatrix}\leq \frac{\I (1+|z|^2)}{(
 \ov{z}-z
)}I_{m_2}.
\end{align}

Compare \eqref{ap14} with the Definition \ref{defWeyl} of the Weyl function in order to
see that $\vp_{\infty}$ is a  Weyl function of \eqref{0.1}, \eqref{0.2}.  Moreover, this Weyl function analytic and contractive in $\BC_-$.
It remains to show that the Weyl function is unique. 

First notice that \eqref{1.3}  yields
\begin{align}\label{ap15}
q(z)W_{k+1}(z)^*jW_{k+1}(z) \geq W_k(z)^*jW_k(z) \quad (k \geq 0).
\end{align}
Thus, we have $q(z)^{k+1}W_{k+1}(z)^*jW_{k+1}(z) \geq j$, and so \eqref{Prop2.2} implies that
\begin{align}\label{ap16}& 
\begin{bmatrix}
 I_{m_1} & 0 
\end{bmatrix}\sum_{k=0}^r
q(z)^kW_k(z)^*C_kW_k(z) \begin{bmatrix}
 I_{m_1} \\ 0 
\end{bmatrix}
\geq (r+1)I_{m_1}.
\end{align}
Therefore, there is an $m_1$-dimensional subspace of vectors $g\in \BC^m$ such that
\begin{align}\label{ap17}& 
\sum_{k=0}^{\infty}
g^*q(z)^kW_k(z)^*C_kW_k(z)g
=\infty .
 \end{align}
The further proof of the uniqueness of the values, which the Weyl function may take
at any fixed $z\in \BC_-$ is easy and coincides  with the arguments in \cite[Theorem 3.8]{FKRS14}.
\end{proof}
 

\begin{flushright}

A.L. Sakhnovich,\\
Fakult\"at f\"ur Mathematik, Universit\"at Wien, \\
Oskar-Morgenstern-Platz 1, A-1090 Vienna, Austria\\
e-mail: oleksandr.sakhnovych@univie.ac.at
\end{flushright}


\begin{thebibliography}{AGKS}

\bibitem{Ci}
 \textit{J.L.~Cieslinski},  
Algebraic construction of the Darboux matrix revisited. --  \textit{J. Phys. A} \textbf{42}  (2009), 404003.

\bibitem{D}
 \textit{P.A.~Deift}, Applications of a commutation formula. --  \textit{Duke Math. J.} \textbf{45} (1978), 267--310.

\bibitem{DFK}
\textit{V.K. Dubovoj, B. Fritzsche, and B. Kirstein}, Matricial version of the classical Schur problem. Teubner-Texte zur
Mathematik [Teubner Texts in Mathematics] \textbf{ 129}, B.G. Teubner Verlagsgesellschaft mbH, Stuttgart, 1992.

\bibitem{FKRS08}
\textit{B. Fritzsche, B. Kirstein, I.Ya. Roitberg, and A.L. Sakhnovich},  Weyl matrix functions and 
inverse problems for discrete Dirac type self-adjoint system: explicit and general 
solutions. --  \textit{Operators and Matrices} \textbf{2} (2008), 201--231.

\bibitem{FKRS14}
\textit{B. Fritzsche, B. Kirstein, I. Ya. Roitberg, and A.L. Sakhnovich}, Discrete Dirac system: rectangular Weyl, functions, direct and inverse problems. -- \textit{Oper. Matrices} \textbf{8}:3 (2014), 799--819.

\bibitem{FKRS17}
\textit{B. Fritzsche, B. Kirstein, I.Ya. Roitberg, and A.L. Sakhnovich},
Stability of the procedure of explicit recovery of skew-selfadjoint Dirac systems from rational Weyl matrix functions.
 -- \textit{Linear Algebra Appl.} \textbf{533} (2017), 428--450.

\bibitem{GeT}
 \textit{F.~Gesztesy and G.~Teschl}, On the double commutation method. --  \textit{Proc. Amer. Math. Soc.} \textbf{124} (1996), 1831--1840.

\bibitem{Gu}
\textit{C. Gu,  H. Hu, and Z. Zhou,} 
Darboux transformations in integrable systems. 
Springer, Dordrecht, 2005.

\bibitem{KS}
 \textit{M.A. Kaashoek and   A.L. Sakhnovich}, Discrete pseudo-canonical system and isotropic Heisenberg magnet. --   \textit{J. Funct. Anal.} 
\textbf{228}  (2005), 207--233.	

\bibitem{KFA}
\textit{R.E.~Kalman, P.~Falb, and M.~Arbib}, Topics in mathematical system theory. 
New York, McGraw-Hill Book Company, 1969.

\bibitem{KoSaTe}
\textit{A. Kostenko, A. Sakhnovich, and G. Teschl},  Commutation methods for Schr\"odinger operators with strongly singular potentials. -- 
\textit{Math. Nachr.} \textbf{285}:4 (2012), 392--410.

 \bibitem{LR}
\textit{P. Lancaster and L. Rodman},  
 Algebraic Riccati equations. Oxford,  Clarendon Press, 1995.

 \bibitem{Mar0}
\textit{V.A. Marchenko},
Stability of the inverse problem of scattering theory. -- 
\textit{Mat. Sb. (N.S.)} \textbf{77(119)}:2 (1968),  139--162. 

\bibitem{Mar}
\textit{V.A. Marchenko}, Nonlinear equations and operator algebras.  D. Reidel, Dordrecht, 1988. 

\bibitem{MS}
 \textit{V.B.~Matveev and M.A.~Salle},  Darboux transformations and solitons.  Springer, Berlin,  1991.

\bibitem{RaRo} 
 \textit{A.C.M. Ran and L. Rodman},  Stability of invariant maximal semidefinite subspaces, II: Applications: selfadjoint rational matrix 
 functions, algebraic Riccati equations. -- \textit{Linear Algebra Appl.} \textbf{63} (1984), 133--173.

 \bibitem{ALS94}
 \textit{A.L.~Sakhnovich}, Dressing procedure for solutions of nonlinear equations and the method of operator identities.
 -- \textit{Inverse Problems} \textbf{10}:3 (1994),  699--710.

\bibitem{ALS16}
\textit{A.L.~Sakhnovich}, Inverse problems for self-adjoint Dirac systems: explicit solutions and stability of the procedure. --  
\textit{Oper. Matrices} \textbf{10}:4 (2016),  997--1008. 

\bibitem{ALSverb}
\textit{A.L.~Sakhnovich}, Verblunsky-type coefficients for Dirac and canonical systems generated by Toeplitz and Hankel matrices, respectively.
arXiv:1711.03064

\bibitem{SaSaR}
 \textit{A.L.~Sakhnovich,  L.A.~Sakhnovich, and I.Ya.~Roitberg},    Inverse problems and nonlinear evolution equations. 
 Solutions, Darboux matrices and Weyl--Titchmarsh functions. De Gruyter Studies in Mathematics \textbf{47},  De Gruyter, Berlin, 2013.

 \bibitem{SaL1}
 \textit{L.A. Sakhnovich}, On  the  factorization  of  the 
transfer matrix
function. -- {Sov. Math. Dokl.} \textbf{17} (1976), 203--207.

 \bibitem{SaL3}
\textit{L.A. Sakhnovich}, Spectral theory of canonical differential
systems, method of operator identities. -- Operator Theory Adv. Appl. {\bf 107},
Birkh\"auser Verlag, Basel, 1999.




 \end{thebibliography}
\end{document}